\documentclass{article}
            \usepackage{graphicx}
            \usepackage{amsmath,amssymb}
            \begin{document}
            \title{Representations of Each Number Type
            that Differ by Scale Factors}
            \author{Paul Benioff\\
            Physics Division, Argonne National
            Laboratory,\\ Argonne, IL 60439, USA\\
            e-mail: pbenioff@anl.gov}

            \maketitle

            \begin{abstract}
            For each type of number, structures that differ by
            arbitrary scaling factors and are isomorphic to one another
            are described. The scaling of number values in one structure,
            relative to the values in another structure, must be
            compensated for by scaling of the basic operations and
            relations (if any) in the structure. The scaling must be
            such that one structure satisfies the relevant number type
            axioms if and only if the other structure does.
            \end{abstract}

            \section{Introduction}
            Numbers play an essential role in many areas of human
            endeavor. Starting with the natural numbers, $N$, of arithmetic,
            one progresses up to integers, $I$, rational numbers, $Ra$, real
            numbers, $R$, and to complex numbers, $C$. In mathematics and
            physics, each of these types of numbers  is referred to as
            \emph{the} natural numbers, \emph{the} integers, rational
            numbers, real numbers, and \emph{the} complex numbers.
            As is well known, though, "the" means "the same up to
            isomorphism" as there are many isomorphic representations
            of each type of number.

            In this paper, properties of different isomorphic
            representations of each number type will be investigated.
            Emphasis is placed on  representations of each number
            type that differ from one another by arbitrary
            scale factors. Here mathematical properties of these
            representations will be described. The possibility that
            these representations for complex numbers may  be
            relevant to physics is described elsewhere \cite{BenNGF}.

            Here the mathematical logical description of a
            representation, as a structure that satisfies a
            set of axioms relevant to the type of system being
            considered \cite{Barwise,Keisler}, is used. For the
            scaled structures considered here, it will be useful
            in some cases to separate the notion of representation
            from that of structure, and consider representations as
            different views of a structure. This will be noted when
            needed.

            Each structure consists of a
            base set, one or more basic operations, basic relations
            (if any), and constants. Any structure  containing a
            base set, basic operations, relations, and constants  that
            are relevant for the number type, and are such that the
            structure satisfies the relevant axioms is a model of
            the axioms. As such it is as good a representation of
            the number type as is any other representation.

            The contents of structures  for the different types of
            numbers and the chosen axiom sets are shown below:
            \begin{itemize}
            \item $\overline{N}=\{N,+,\times,<, 0,1\}$\hspace{1cm}
            \parbox[t]{6.5cm}{Nonnegative elements of a discrete
            ordered commutative ring with identity
            \cite{Kaye}.}

            \item $\overline{I}=\{I,+,-,\times,<, 0,1\}$\hspace{1cm}
            Ordered integral domain \cite{integer}.

            \item $\overline{Ra}=\{Ra,+,-,\times,\div,<, 0,1\}$\hspace{0.8cm}
            Smallest ordered field \cite{rational}.

            \item $\overline{R}=\{R,+,-,\times,\div,<, 0,1\}$\hspace{1cm}
            Complete ordered field \cite{real}.

            \item $\overline{C}=\{C,+,-,\times,\div,^{*}, 0,1\}$\hspace{0.5cm}
            \parbox[t]{6.5cm}{Algebraically closed field of characteristic $0$ plus
            axioms for complex conjugation \cite{complex,comcon}.}
            \end{itemize}

            Here an overline, such as in $\overline{N},$ denotes a
            structure. No overline, as for $N$, denotes a base set.
            The complex conjugation operation has been added as a
            basic operation to $\overline{C}$ as it makes the
            development much easier. The same holds for the
            inclusion of the division operation in
            $\overline{Ra},$ $\overline{R},$ and $\overline{C}.$

             For this work, the choice of  which axioms are used for
            each of the number types is not important. For example,
            an alternate choice for $\overline{N}$ is to use the axioms
            of arithmetic \cite{Smullyan}. In this case $\overline{N}$
            is changed by deleting the constant $1$ and adding a
            successor operation.  There are also other axiom
            choices for the real numbers \cite{Lang}.

            The importance of the axioms is that they will be used
            to show that, for two structures related by a scale factor,
            one satisfies the axioms if and only if the other does.
            This is equivalent to showing that one is a structure for
            a given number type if and only if the other one is a
            structure for the same number type.

            These ideas will be expanded in the following sections.
            The next section gives a general treatment of fields. This
            applies to all the number types that satisfy the field
            axioms (rational, real, complex numbers).
            However much of the section applies to other
            numbers also (natural numbers, integers). The following
            five sections apply the general results to each of the
            number types. The discussions are mainly limited to
            properties of the number type that are not included in
            the description of fields.

            Section \ref{NTSC}  expands the descriptions of the
            previous sections by considering
            $\overline{N},\overline{I},\overline{Ra},\overline{R}$
            as substructures of $\overline{C}.$ In this case the
            scaling factors relating two structures of the same
            type are complex numbers.

            Section \ref{D} concludes the paper with a discussion of
            some consequences and possible uses of  these
            representations in physics.

            \section{General Description of Fields}\label{GD}
            It is useful to describe the results of this work for
            fields in general. The results can then be applied to
            the different number types, even those that are not fields.
            Let $\overline{S}$ be a field structure where
             \begin{equation}\label{S}\overline{S}=\{S,+,-,\times,
            \div,0,1\}.\end{equation}Here $S$  with no overline denotes a base set,
            $+,-,\times,\div$ denote the basic field operations, and $0,1$
            denote constants. Denoting $\overline{S}$ as a field
            structure implies that $\overline{S}$ is a structure that
            satisfies the axioms for a field \cite{field}.

            Let $\overline{S}_{p}$ where \begin{equation}\label{Sp}\overline{S}_{p}
            =\{S,+_{p},-_{p},\times_{p},\div_{p},0_{p},1_{p}\}.\end{equation} be
            another structure on the same set $S$ that is in $\overline{S}.$
            The idea is to require that $\overline{S}_{p}$ is also a field structure
            on $S$ where the field values of the elements of $S$
            in $\overline{S}_{p}$ are scaled
            by $p,$ relative to the  field values in $\overline{S}.$ Here $p$
            is a field value in $\overline{S}.$

            The goal is to show that this is possible in that one can define
            $\overline{S}_{p}$ so that $\overline{S}_{p}$ satisfies the field
            axioms if and only if $\overline{S}$ does. To this end the notion of
            correspondence is introduced as a relation between the field values
            of $\overline{S}_{p}$ and $\overline{S}.$ The field value, $a_{p},$
            in $\overline{S}_{p}$ is said to \emph{correspond} to the field value, $pa,$
            in $\overline{S}.$  As an example, the identity value, $1_{p},$ in
            $\overline{S}_{p}$ corresponds to the value $p\times 1=p$ in $\overline{S}.$

            This shows that correspondence is distinct from the concept of
            sameness.  $a_{p}$ is the same value in $\overline{S}_{p}$
            as  $a$ is in $\overline{S}.$ This differs from $pa$ by the
            factor $p.$ The distinction between correspondence and sameness
            is present only if $p\neq 1.$  If $p=1$, then the two concepts
            coincide, and $\overline{S}_{p}$ and $\overline{S}$
            are the same structures as far as scaling is concerned.

            So far a scaling factor has been introduced that relates field
            values between $\overline{S}_{p}$ and $\overline{S}.$ This must
            be compensated for by a scaling of the basic operations in
            $\overline{S}_{p}$ relative to those of $\overline{S}.$

            The correspondences of the basic field operations and
            values in $\overline{S}_{p}$ to those in $\overline{S}$ are given by,
            \begin{equation}\label{opsp1}\begin{array}{c}a_{p}=pa\\+_{p}=+,
            \hspace{1cm}-_{p}=-\\\times_{p}=\frac{\textstyle\times}{\textstyle
            p},\hspace{1cm}\div_{p}=p\div.\end{array}\end{equation}

            One can use these scalings to replace the basic operations and
            constants in $\overline{S}_{p}$ and define $\overline{S}^{p}$ by,
            \begin{equation}\label{Sp11}\overline{S}^{p}=\{S,+,-,\frac{\times}{p},
            p\div,0,p\}.\end{equation} Here the subscript, $p,$
            in $\overline{S}_{p},$ Eq. \ref{Sp} is replaced by $p$
            as a superscript to distinguish $\overline{S}^{p}$
            from $\overline{S}_{p}.$

            Both $\overline{S}^{p}$ and $\overline{S}_{p}$ can
            be considered  as different representations or views
            of a structure that differs from $\overline{S}$ by a
            scaling factor, $p$. A useful
            expression of the relation between $\overline{S}_{p}$
            and $\overline{S}^{p}$ is that $\overline{S}^{p}$
            is  referred to either as a representation
            of $\overline{S}_{p}$ on $\overline{S},$ or as an explicit
            representation of $\overline{S}_{p}$ in terms of the
            operations and element values of $\overline{S}.$

            Besides changes in the definitions of the basic operations
            given in Eq. \ref{opsp1} and distinguishing between
            correspondence and sameness, scaling introduces another change.
            This is that one must drop the usual assumption that the
            elements of the base set, $S,$ have fixed values, independent
            of structure membership. Here the field values of the elements
            of $S$, with one exception, depend on the structure containing
            $S.$ In particular, to say that
            $a_{p}$ in $\overline{S}_{p}$ corresponds to $pa$ in $\overline{S}$
            means that the element of $S$ that has the value $a_{p}$ in
            $\overline{S}_{p}$ has the value $pa$ in $\overline {S}.$ This is
            different from the element of $S$ that has the same value, $a,$
            in $\overline{S}$ as $a_{p}$ is in $\overline{S}_{p}.$

\begin{figure}[h]\begin{center}
           \resizebox{100pt}{100pt}{\includegraphics[250pt,200pt]
           [500pt,500pt]{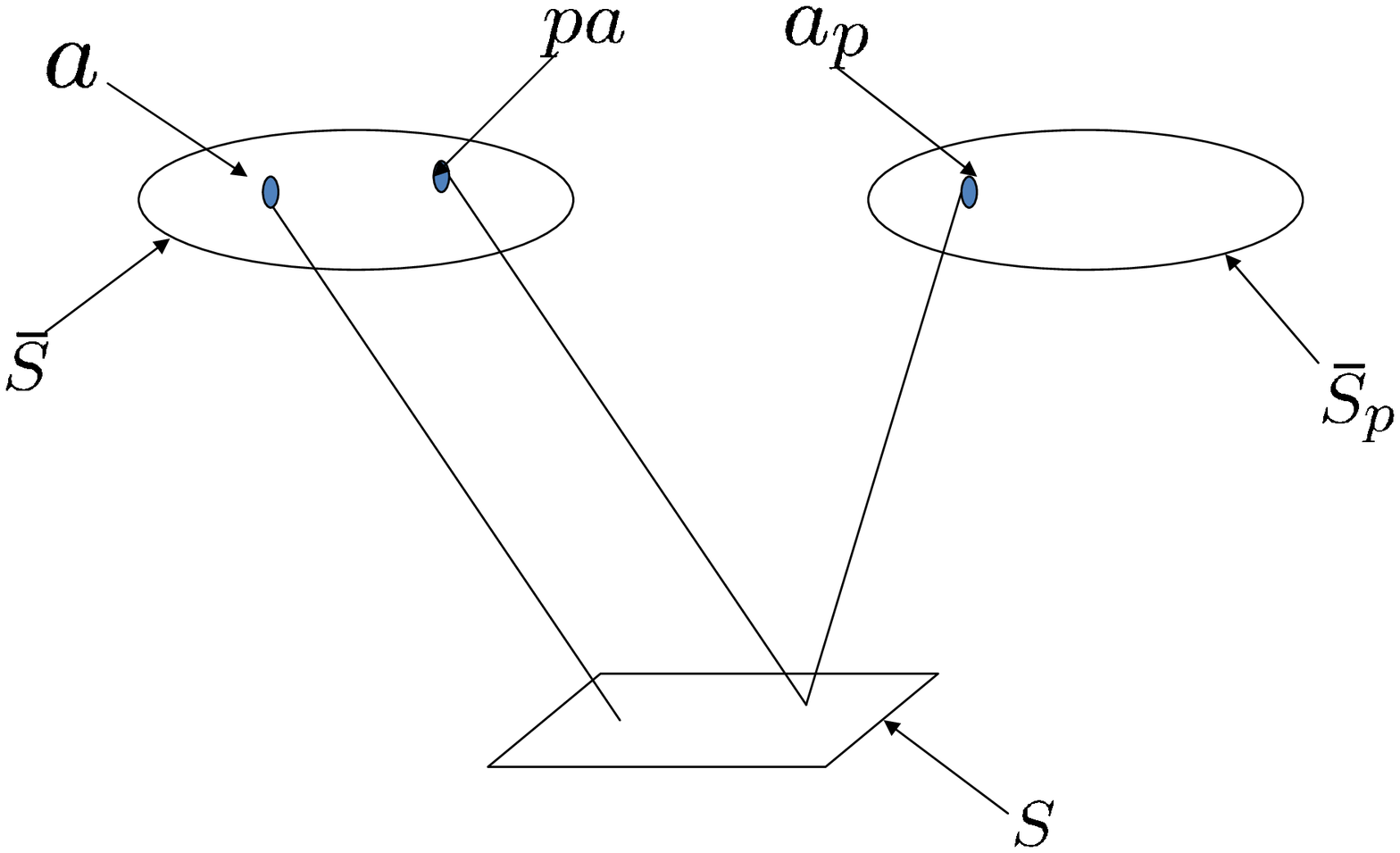}}\end{center}
           \caption{Relations between Elements in the base set $S$
           and their  Values in  the Structures $\bar{S}$
           and $\bar{S}_{p}.$ Here $a_{p}$ is the same value in
           $\bar{S}_{p}$ as $a$ is in $\bar{S}.$ The lines show that
           they are values for different elements of $S.$ The
           lines also show that the $S$ element that has the value $a,$ as
           $a_{p}$ in $\bar{S}_{p}$, has the value $pa$
           in $\bar{S}.$}\label{RENT1}\end{figure}

            These relations are shown schematically in Figure
            \ref{RENT1}. The valuations associated with elements in
            the base set $S$ are shown by lines  from $S$ to the
            structures $\bar{S}$ and $\bar{S}_{p}.$

            The one exception to the structure dependence of valuations
            is the element of $S$ with the value $0.$  This value-element
            association is fixed and is independent of all values of $p.$
            In this sense it is the "number vacuum" as it is invariant
            under all changes $p\rightarrow p'.$\footnote{Like the physical
            vacuum which is unchanged under all space time translations.}

            Another approach to understanding the differences
            between $\overline{S},\overline{S}^{p},$ and
            $\overline{S}_{p},$ is to distinguish carefully how
            structure elements are described, when viewed from inside
            and from outside the structure.\footnote{The importance of
            distinguishing between inside and outside views of structures
            is well known in mathematical logic as it plays a role in
            the resolution of the  Skolem paradox \cite{LS}.}
            Comparison of different structures necessarily is done
            from outside the structures.

            Here there are two structures, $A,B$, with $A$ represented by
            $\overline{S}$ and $B$ by both $\overline{S}^{p}$
            and $\overline{S}_{p}$. There is just one representation,
            $\overline{S},$ for $A$ because the inside and outside
            views coincide for this structure. For $A$, correspondence
            and sameness coincide. In mathematical logical terms the
            properties of the elements of $A$ are \emph{absolute}
            \cite{absolute}. They are the same when viewed inside
            or outside the structure. For this reason
            $\overline{S}$ will often be referred to as a
            structure instead of a representation.

            The situation is quite different for the structure, $B$.
            Here $\overline{S}^{p}$ and $\overline{S}_{P}$ represent,
            respectively, external and internal views of $B$. For
            example, the element of $S$ that has value $p$ in
            $\overline{S},$ also has value $p$ in $B$ when viewed
            externally.  This shown explicitly in the representation
            $\overline{S}^{p}.$

            When viewed from inside $B$, the element of $S$ that has
            value $p$ when viewed externally, has the value $1$ when
            viewed internally or inside $B$.  This is designated by
            the identity $1_{p}$ in $\overline{S}_{p}.$

            The situation is somewhat different for the basic
            operations of multiplication and division. The basic
            operation of multiplication in $B,$ is shown internally
            by $\times_{p}$ in $\overline{S}_{p}.$ It is seen externally
            as the operation $\times/p.$  Here $\times/p=
            \times(-)\div p$ is  the external and internal view
            of an operation in $A.$ As seen in $\overline{S},$
            $\times$ and $\div$ denote the multiplication and
            division operations in $A.$

            One can also define isomorphic maps between the
            structure representations.  Define the maps $W^{p}$ and
            $W_{p}$ by \begin{equation}\label{WpWp}\overline{S}_{p}
            =W_{p}\overline{S}^{p}=W_{p}W^{p}\overline{S}=F_{p}
            \overline{S}.\end{equation}$W^{p}$ maps
            $\overline{S}$ onto $\overline{S}^{p}$ and $W_{p}$
            maps $\overline{S}^{p}$ onto $\overline{S}_{p}.$ $W^{p}$
            is a map from one structure to another, and $W_{p}$ is a
            map between different representations of the same structure.

            $W^{p}$ and $W_{p}$ are defined by\begin{equation}\label{Wp}
            \begin{array}{c}W^{p}(a)=pa,\\ W^{p}(a\pm b)= W^{p}(a)W^{p}(\pm )W^{p}(b)
            =(pa)\pm (pb)\\W^{p}(a\times b)=W^{p}(a)W^{p}(\times)W^{p}(b)
            =(pa)\frac{\textstyle \times}{\textstyle p}(pb)\\W^{p}(a\div b)=W^{p}(a)
            W^{p}(\div)W^{p}(b)=(pa)(p\div)(pb)\end{array}\end{equation}
            and\begin{equation}\label{Wpp}\begin{array}{c}W_{p}(pa)=a_{p},
            \\W_{p}(pa\pm pb)=W_{p}(pa)W_{p}(\pm )W_{p}(pb)
            =a_{p}\pm_{p}b_{p}\\W^{p}(pa\frac{\textstyle
            \times}{\textstyle p}pb)= W^{p}(pa)W^{p}(\frac{\textstyle \times}
            {\textstyle p})W^{p}(pb) =a_{p}\times_{p}b_{p}\\W^{p}
            (pa(p\div)pb)=W^{p}(pa)W^{p}(p\div)W^{p}(pb)=a_{p}\div_{p}b_{p}.
            \end{array}\end{equation}

            It is important to emphasize that the definition of
            $\overline{S}^{p}$ is not just a relabeling of the
            elements of $\overline{S}.$ One way to see this is
            to show that the description of the relations between
            $\overline{S}^{p}$ and $\overline{S}$ by use of
            isomorphisms is necessary but not sufficient.  To see
            this let $\overline{S}^{wyz}$ denote a structure where
            \begin{equation}\label{Swyz}\overline{S}^{wyz}=\{+,-,
            \frac{\times}{w},(y\div),0,z\}.\end{equation} Here
            $w,y,z$ are arbitrary field values in $\overline{S}.$

            Define the isomorphism $W^{wyz}:\overline{S}\rightarrow
            \overline{S}^{wyz}$ by\begin{equation}\label{Wwyzdef}
            \begin{array}{c}W^{wyz}a=za,\hspace{.5cm}W^{wyz}(\pm)
            =\pm\\W^{wyz}(\times)=\frac{\textstyle\times}{\textstyle w},
            \hspace{.5cm}W^{wyz}(\div)=y\div.\end{array}\end{equation}

            It is clear from the definition of isomorphisms that
            $\overline{S}^{wyz}$ satisfies the field axioms if and
            only if $\overline{S}$ does. This follows from the
            observation that all the field axioms \cite{field} are
            equations. For example, the existence of a multiplicative
            identity, $a\times 1 = a$ in $\overline{S}$ gives the
            equivalences\begin{equation}\label{wyzeqs}\begin{array}{l}
            a\times 1=a\Leftrightarrow W^{wyz}(a)W^{wyz}(\times)
            W^{wyz}(1)=W^{wyz}(a)\\\hspace{1cm}\Leftrightarrow
            za\frac{\textstyle \times}{\textstyle w}(z) =za.\end{array}
            \end{equation}The righthand equation shows that $z$ is
            the multiplicative identity in $\overline{S}^{wyz}.$

            The equivalences of Eq. \ref{wyzeqs} show that if one
            attempts to interpret $\overline{S}^{wyz}$ as an external
            view of a different structure whose operations and constants
            are defined in terms of those in $\overline{S}$, then
            $\overline{S}^{wyz}$  does not represent a field. Requiring
            that $\overline{S}^{wyz}$ represents a field, limits one to
            concluding that  $\overline{S}^{wyz}$ and $\overline{S}$
            represent the same structure.  In particular, $\overline{S}^{wyz}$
            is just a relabeling of the elements of $\overline{S}$ with
            no different valuations implied. Note that this
            limitation is not present for the case where $w=y=z.$

            The following  two theorems summarize the relations
            between $\overline{S}^{p},\overline{S}_{p},$ and
            $\overline{S}.$. The first theorem shows the
            invariance of  equations under the maps, $W^{p}$ and
            $W_{p}$ where $p,$ is a scaling factor. It also shows,
            that the correspondence between between element values
            in $\overline{S}_{p}$ and those in $\overline{S},$
            extends to general terms.
            \newtheorem{Sax}{Theorem}\begin{Sax}\label{Sterms}
            Let $t$ and $u$ be terms in $\overline{S}.$ Let
            $\overline{S}_{p}$  and $\overline{S}^{p}$ be
            as defined in Eqs. \ref{Sp} and \ref{Sp11}.
            Then $t=u\Leftrightarrow t^{p}=u^{p}\Leftrightarrow t_{p}=u_{p}$
            where $t^{p}=W^{p}t,$ $u^{p}= W^{p}u,$ and $W_{p}t^{p}=
            t_{p},W_{p}u^{p}=u_{p}.$\end{Sax}\underline{Proof}: It follows
            from the properties of $W^{p}$ and $W_{p},$ Eqs. \ref{Wp}
            and \ref{Wpp}, that $t^{p}=W^{p}t=pt$ and
            $u^{p}=W^{p}u=pu.$ Also $t_{p}=W_{p}t^{p}$ and
            $u_{p}=W_{p}u^{p}.$ This gives\begin{equation}\label{tpup}
            t_{p}=u_{p}\Leftrightarrow t^{p}=u^{p}\Leftrightarrow pt=pu
            \Leftrightarrow t= u.\end{equation}From the left the first
            equation is in $\overline{S}_{p}$, the second in
            $\overline{S}^{p},$ the third in both $\overline{S}^{p}$
            and $\overline{S},$ and the fourth in $\overline{S}.$

             It remains to see in detail that the correspondence between
             element values in $\overline{S}_{p}$  and those in
             $\overline{S}$ extends to terms. Let\begin{equation}
            \label{termpx} t_{p}=(\sum_{j,k=1}^{m})_{p}\frac{(a_{p})^{j}}
            {(b_{p})^{k}}\mbox{}_{p}.\end{equation} The external view
            of $t_{p}$ in $\overline{S}^{p}$ is \begin{equation}
            \label{termpxx} t^{p}=(\sum_{j,k=1}^{m})^{p}\frac{(pa)^{j}}
            {(pb)^{k}}\mbox{}^{p}.\end{equation}In
            the numerator, the $j$ $pa$ factors and $j-1$
            multiplications contribute factors of $p^{j}$ and $p^{-j+1},$
            respectively, to give a factor $p$. This is canceled by
            a similar factor arising from the denominator.
            A final factor of $p$ arises from the representation
            of the solidus as shown in Eq. \ref{opsp1} for division.

            Using this and the fact that addition is not scaled,
            gives the result that \begin{equation}\label{termpxxtx}
            t^{p}=(\sum_{j,k=1}^{m})^{p}\frac{(pa)^{j}}{(pb)^{k}}\mbox{}^{p}
            =p(\sum_{j,k=1}^{m})\frac{(a)^{j}}{(b)^{k}}\mbox{}=pt.
            \end{equation} Eq.\ref{tpup} and the theorem follow
            from the fact that Eq. \ref{termpxxtx} holds for any
            term, including $u^{p}.$ $\blacksquare$

            From this one has\newtheorem{SAx}[Sax]{Theorem}\label{S}\begin{SAx}
            $\overline{S}_{p}$ satisfies the field axioms if and only
            if $\overline{S}^{p}$ satisfies the field axioms if and
            only if $\overline{S}$ satisfies the field axioms.\end{SAx}
            \underline{Proof}: The theorem follows from Theorem 1 and
            the fact that all the field axioms are equations for terms.
            $\blacksquare$

            The constructions described here can be iterated. Let $p$
            be a number value in $\overline{S}$ and $q\equiv q_{p}$
            be a number value in $\overline{S}_{p}.$ Let $\overline{S}_{q|p}$
            be another field structure on the base set, $S,$ and let
            $\overline{S}^{q|p}$ be the representation of
            $\overline{S}_{q|p}$ using the operations and constants
            of $\overline{S}_{p}.$\footnote{An equivalent way to
            define $\overline{S}^{q|p}$ is as the representation
            of $\overline{S}_{q|p}$ on $\overline{S}_{p}.$} In more
            detail \begin{equation}\label{Sqpdown}\overline{S}_{q|p}
            =\{S,\pm_{q|p},\times_{q|p},\div_{q|p},0,1_{q|p}\}
            \end{equation} and \begin{equation}\label{Sqpup}
            \overline{S}^{q|p}=\{S,\pm_{p}, \frac{\times_{p}}{q},q\div_{p},
            0,q1_{p}\}.\end{equation}$\overline{S}_{q|p}$ is related to
            $\overline{S}_{p}$ by the scaling factor $q.$ The goal is to
            determine the scaling factor for the representation of
            $\overline{S}_{q|p}$ on $\overline{S}.$

            To determine this, let  $a_{q|p}$ be a value
            in $\overline{S}_{q|p}.$ This corresponds to a value
            $q_{p}\times_{p}a_{p}$ in $\overline{S}_{p}.$ Here
            $a_{q|p}$ is the same value in $\overline{S}_{q|p}$
            as $q_{p}\times_{p}a_{p}$ is in $\overline{S}^{q|p}$
            as $a_{p}$ is in $\overline{S}_{p}.$

            The value in $\overline{S}$ that corresponds to $a_{q|p}$
            in $\overline{S}_{q|p}$ can be determined from its correspondent,
            $q_{p}\times_{p}a_{p},$ in $\overline{S}_{p}.$ The value in
            $\overline{S}$ that corresponds
            to  $q_{p}\times_{p}a_{p}$ is obtained by use of
            Eqs. \ref{opsp1} and \ref{Wpp}. It is given by
            \begin{equation}\label{Wypinv}(W_{p})^{-1}(q_{p}\times_{p}a_{p})=
            (pq)\frac{\textstyle \times}{\textstyle p}(p
            a)=pqa.\end{equation}Here $q$ and $a$
            are the same values in $\overline{S}$  as $q_{p}$ and $a_{p}$
            are in $\overline{S}_{p}.$ Also $(W_{p})^{-1},$ as the inverse
            of $W_{p},$ maps $\overline{S}_{p}$ onto $\overline{S}^{p}.$

            This is the desired result because it shows that
            two steps, first with $p$ and then with $q$ is
            equivalent to one step with $qp.$  This result shows that
            the representation of $\overline{S}_{q|p}$ on
            $\overline{S}$ is given by\begin{equation}\label{SqpS}
            \overline{S}^{q|p}=\{S,\pm, \frac{\times}{qp},qp\div,
            0,qp1\}.\end{equation}

            Note that the steps commute in that the same result
            is obtained if one scales first by $q$ and then by $p$
            as $pq =qp.$  Here $qp$ is a value in $\overline{S}.$
            Also this is equivalent to determining
            the scale factor for $\overline{S}_{q}$ on $\overline{S}$
            provided one accounts for the fact that $q$ is a value in
            $\overline{S}_{p}$ and not in $\overline{S}.$

            \section{Natural Numbers}\label{NN} The natural numbers
            differ from the generic representation in that they are
            not fields \cite{Kaye}.  This is shown by the structure
            representation for $\overline{N},$\begin{equation}
            \label{barN}\overline{N}=\{N,+,\times,<,0,1\}.\end{equation}

            The structure corresponding to $\overline{S}_{n}$ is
            \begin{equation}\label{Nn}\overline{N}_{n}
            =\{N_{n},+_{n},\times_{n},<_{n},0_{n},1_{n}\}.
            \end{equation} Here $n$ is any natural number $>0.$

             One can use Eq. \ref{opsp1} to represent
            $\overline{N}_{n}$ in terms of the basic operations,
            relations, and constants of $\overline{N}.$  It is
            \begin{equation}\label{Nn11}\overline{N}^{n}=
            \{N_{n},+,\frac{\times}{n},<,0,n\}.
            \end{equation} Note that, as is the case for addition,
            the order relation is the same in $\overline{N}_{n}$
            as in $\overline{N}^{n}$ and in $\overline{N}.$ As was
            the case for fields, $\overline{N}^{n}$ and
            $\overline{N}_{n}$ represent external and internal views
            of a different structure than that represented by
            $\overline{N}$.  For $\overline{N}$ the external and
            internal views coincide.

            The multiplication operator in $\overline{N}^{n},$
            $\times/n$  has the requisite properties. This
            can be seen by the equivalences between multiplication
            in $\overline{N}_{n},\overline{N}^{n},$ and
            $\overline{N}:$  $$a_{n}\times_{n}b_{n}
            =c_{n}\Leftrightarrow na\frac{\times}{n}n
            b=nc\Leftrightarrow a\times b=c.$$

            Note that the simple verification of these equivalences takes
            place outside the three structures and not within any natural
            number structure.  For this reason division by $n$ can be used
            to verify the equivalences even though it is not part of any
            natural number structure.

            These equivalences show that $n$ is the multiplicative
            identity in $\overline{N}^{n}$ if and only if $1$
            is the multiplicative identity in $\overline{N}.$ To
            see this set $b=1.$

            The structure, $\overline{N}^{n},$  Eq. \ref{Nn11}, and that of Eq.
            \ref{Nn}, differ from  the generic description, Section \ref{GD},
            in that the base set $N_{n}$ is a subset of $N.$ $N_{n}$ contains
            just those elements of $N$ whose values in $\overline{N}$ are
            multiples of $n.$ For example, the element with value $n$
            in $\overline{N}$ has value $1$ in $\overline{N}_{n},$ and
            the element with value $na$ in $\overline{N}$ has value $a$
            in $\overline{N}_{n}.$ Elements with values $na+l$ in
            $\overline{N}$ where $0<l<n,$ are absent from $N_{n}.$\footnote{The
            exclusion of elements of the base set in different
            representations occurs only for the natural numbers and
            the integers. It is a consequence of their not being
            closed under division.}

            As noted before, the choice of the representations of the
            basic operations and relation and constants in
            $\overline{N}_{n},$ as shown in $\overline{N}^{n},$ is
            determined by the requirement that $\overline{N}^{n}$
            satisfies the natural number axioms \cite{Kaye} if and only if
            $\overline{N}$ satisfies the axioms. For the axioms  that are equations
            and do not use the ordering relation, this follows immediately from
            Theorem 1.

            For axioms that use the order relation, the requirement
            follows from the fact that $<_{n}=<$ and  for any
            pair of terms $t^{n},u^{n}$,
            $$t^{n}<_{n}u^{n}\Leftrightarrow nt<nu\Leftrightarrow
            t<u.$$ Here Eq. \ref{termpxxtx} was used with
            $b=1.$ Note that $n>0.$

            It follows from these considerations that
            \newtheorem{Nax}[Sax]{Theorem}\label{N}\begin{Nax}$\overline{N}_{n}$
            satisfies the axioms of arithmetic if and only if
            $\overline{N}^{n}$ does if and only if
            $\overline{N}$ does.\end{Nax}. A simple example that
            illustrates the theorem is the axiom of discreteness for
            the ordering, $0<1\bigwedge\forall a(a>0\rightarrow
            a\geq 1)$ \cite{Kaye}. One has the equivalences $$\begin{array}{l}
            0<_{n}1_{n}\bigwedge\forall a_{n}(a_{n}>_{n}0\rightarrow
            a_{n}\geq_{n} 1_{n})\Leftrightarrow 0<n\bigwedge\forall
            a(na>0\rightarrow na\geq n1)\\\hspace{1cm}\Leftrightarrow
            0<1\bigwedge\forall a(a>0\rightarrow a\geq 1).\end{array}$$ Subscripts
            are missing on $0$ because the value remains the same in all structures.

            \section{Integers}\label{I}Integers generalize the
            natural numbers in that negative numbers are included.
            Axiomatically they can be characterized as an ordered
            integral domain \cite{integer}.  As a structure,
            $\overline{I}$ is given by\begin{equation}\label{I}
            \overline{I}=\{I,+,-,\times,<,0,1\}.
            \end{equation} $I$ is a base set, $+,-,\times$ are the
            basic operations, $<$ is an order relation, and $0,1$
            are the additive and multiplicative identities.

            Let $j$ be a positive integer. Let $\overline{I}_{j}$ be the
            structure \begin{equation}\label{Ij}\overline{I}_{j}=\{I_{j},
            +_{j},-_{j},\times_{j},<_{j},0_{j},1_{j}\}.\end{equation}$I_{j}$
            is the subset of $I$ containing all and only those
            elements of $I$ whose values in $\overline{I}$ are positive
            or negative multiples of $j$ or $0.$

            The representation of $\overline{I}_{j}$ in terms of elements,
            operations, and relations in $\overline{I}$ is given by
            \begin{equation}\label{IjI}\overline{I}^{j}=\{I_{j},
            +,-,\frac{\times}{j},<,0,j\}.\end{equation}This structure
            differs from that of the natural numbers, Eq. \ref{Nn11},
            by the presence of the additive inverse, $-.$

            The proof that $\overline{I}_{j}$ and $\overline{I}^{j}$
            satisfy the integer axioms if and only if $\overline{I}$
            does, is similar to that for Theorem \ref{N}. The only
            new operation is the additive inverse. Since the axioms
            for this are similar to those already present, details of the
            proof for axioms involving subtraction will be skipped.

            A new feature enters in the case that $j$ is negative.
            It is sufficient to consider the case where $j=-1$ as
            the case for other negative integers can be described as
            a combination of $j=-1$ followed by  scaling with a
            positive $j$. This would be done by extending the
            iteration process, described for fields in Section \ref{GD},
            to integers.

            The integer structure representations for $j=-1$
            that correspond to Eqs. \ref{Ij} and \ref{IjI} are
            given by \begin{equation}\label{I-1}\overline{I}_{-1}=
            \{I,+_{-1},-_{-1},\times_{-1},<_{-1},0_{-1},1_{-1}\},
            \end{equation} and \begin{equation}\label{I-1I}\overline{I}^{-1}
            =\{I,+,-,\frac{\times}{-1},>,0,-1\}.\end{equation}
            The main thing to note here is that the order
            relation, $<_{-1},$ in $\overline{I}_{-1}$ corresponds
            to $>$ in $\overline{I}.$

            $\overline{I}^{-1}$ can also be described as a
            reflection of  the whole structure, $\overline{I}=
            \overline{I}_{1}$ through the origin at $0.$
            Not only are the integer values reflected but also the
            basic operations and order relation are reflected.
            Also in this case the base set, $I_{-1}$ =$I.$

            Eq. \ref{I-1I} indicates that $-1$ is the identity
            and $-1$ is positive in $\overline{I}^{-1}.$ These follow from
            $$a_{-1}\times_{-1}1_{-1}=a_{-1}\Leftrightarrow
            (-a)\frac{\times}{-1}(-1)=-a$$ and$$0<_{-1}1_{-1}\Leftrightarrow 0>-1
            .$$ This equivalence shows that the relation, $>$, which is
            interpreted as greater than in $\overline{I},$ is interpreted
            as less than in $\overline{I}^{-1}.$ Thus, as a relation in
            $\overline{I}^{-1},$ $0>-1$ means $-1$ is greater than $0.$
            As a relation in $\overline{I}$, $0>-1$ means $-1$ is less than $0.$

            This is an illustration of the relation of the
            ordering relation $<_{-1}$ to $<.$ Integers which
            are positive in $\overline{I}_{-1}$ and $\overline{I}^{-1},$
            are negative in $\overline{I}.$ It follows that
            $0<_{-1}1_{-1}<_{-1}2_{-1}<_{-1}\cdots$ is
            true in $\overline{I}_{-1}$ if and only if $0>-1>-2>\cdots$
            is true in $\overline{I}^{-1}$ and in $\overline{I}.$

            These considerations show that $\overline{I}_{-1}$
            and $\overline{I}^{-1}$ satisfy the integer
            axioms if and only if $\overline{I}$ satisfies the axioms.
            For axioms not involving the order relation the proof is
            similar to that for $\overline{I}_{j}$ for $j>0.$ For axioms
            involving the order relation  the proofs proceed by restating
            axioms for $\overline{I}$ in terms of $>.$

            An example is the axiom for transitivity $a<b\bigwedge
            b<c\Rightarrow a<c.$ For this axiom the validity of
            the equivalence $$\begin{array}{l}(a_{-1}<_{-1}b_{-1}
            \bigwedge b_{-1}<_{-1}c_{-1}\Rightarrow a_{-1}<_{-1}c_{-1})
            \\\hspace{1cm}\Leftrightarrow (-a>-b\bigwedge -b>-c
            \Rightarrow -a>-c)\end{array}$$ shows that this axiom
            is true in $\overline{I}_{-1}$  if and only if it is
            true in $\overline{I}^{-1}$ if and only if it is true
            in $\overline{I}.$

            These considerations are sufficient to show that a
            theorem equivalent to Theorem \ref{N} holds for
            integers:\newtheorem{Iax}[Sax]{Theorem}\label{I}
            \begin{Iax}For any integer $j\neq 0$, $\overline{I}_{j}$
            satisfies the integer axioms if and only if
            $\overline{I}^{j}$ does if and only if
            $\overline{I}$ does.\end{Iax}

            \section{Rational Numbers}\label{RaN} The next type of
            number to consider is that of the rational numbers.
            Let $\overline{Ra}$ denote a rational number structure
            \begin{equation}\label{Ra}\overline{Ra}=\{Ra,+,-,\times,
            \div,<, 0,1\}.\end{equation} For each positive rational number
            $r$ let $\overline{Ra}_{r}$ denote the structure\begin{equation}
            \label{Rar}\overline{Ra}_{r}=\{Ra,+_{r},-_{r},\times_{r},
            \div_{r},<_{r}, 0_{r},1_{r}\}.\end{equation}

            Eqs. \ref{Ra} and \ref{Rar} show that
            $\overline{Ra}$ and $\overline{Ra}_{r}$ have the same
            base set, $Ra.$ This is a consequence of the fact that
            rational numbers are a field. As such $\overline{Ra}$
            nd  $\overline{Ra}_{r}$ are special
            cases of the generic fields described in Section
            \ref{GD}. Note also that $\overline{Ra}$ is the same as
            $\overline{Ra}_{1}.$

            The definition of $\overline{Ra}_{r}$ is made specific
            by the representation of its elements in terms of those
            of $\overline{Ra}.$ It is\begin{equation}\label{RarRa}
            \overline{Ra}^{r}=\{Ra,+,-,\frac{\times}{r}, r\div,<,
            0,r\}.\end{equation}

            As was seen for fields in Section
            \ref{GD}, the number values of the elements of $Ra$
            depend on the structure containing $Ra.$ The element of
            $Ra$ that has value $ra$ in $\overline{Ra}$ has value
            $a$ in $\overline{Ra}^{r}.$ The element of $Ra$ that
            has the value $a$ in $\overline{Ra}^{r}$ is different from
            the element that has the same value $a$ in $\overline{Ra}.$
            The only exception is the element with value $0$ as this
            value is the same for all $\overline{Ra}_{r}.$ Also, as
            noted in Section \ref{GD}, $\overline{Ra}^{r}$ and
            $\overline{Ra}_{r}$ represent external and internal
            views of a structure that differs from that represented
            by $\overline{Ra}.$

            As was the case for multiplication, the relation between
            $\div_{r}=r\div$ and $\div$ is fixed by the requirement
            that $\overline{Ra}_{r}$ satisfy the rational number
            axioms \cite{rational} if and only if $\overline{Ra}^{r}$
            satisfies the axioms if and only if $\overline{Ra}$ does.
            this can be expressed as a theorem:
            \newtheorem{raAx}[Sax]{Theorem}
            \begin{raAx} Let $r$ be any nonzero rational number.
            Then $\overline{Ra}$ satisfies the rational number axioms
            if and only if $\overline{Ra}^{r}$ satisfies the axioms
            if and only if $\overline{Ra}_{r}$ does.\end{raAx}
            \underline{Proof}:\\Since the axioms for rational
            numbers include those of an ordered field, the proof
            contains a combination of that already given for fields
            in Section \ref{GD}, Theorem \ref{S}, and for the ordering
            axioms for integers, as in Theorem \ref{I}. As a result
            it will not be repeated here.

            It remains to prove that $\overline{Ra}_{r}$ is the
            smallest ordered field if and only if $\overline{Ra}^{r}$
            is if and only if $\overline{Ra}$ is. To show this,
            one uses the isomorphisms defined in Section \ref{GD}.

            Let $\overline{S}$ be an ordered field.  Let
            $W^{r}$ and $W_{r}$ be isomorphisms whose definitions
            on $\overline{Ra}$ and $\overline{Ra}^{r}$ follow
            that of $W^{p}$ and $W_{p}$ in Eqs. \ref{Wp} and \ref{Wpp}.
            That is\begin{equation}\label{RarWr}\overline{Ra}_{r}=
            W_{r}\overline{Ra}^{r}=W_{r}W^{r}\overline{Ra}.\end{equation}

            Since these maps, as isomorphisms, are one-one onto
            and are order preserving, they have inverses which
            are also isomorphisms. In this case one has
            \begin{equation}\label{RaS}\begin{array}{l}
            \overline{Ra}\subseteq\overline{S}\Rightarrow
            (W^{r})\overline{Ra}=\overline{Ra}^{r}
            \subseteq\overline{S}\\\hspace{1cm}\Rightarrow
            (W_{r})\overline{Ra}^{r}=\overline{Ra}_{r}
            \subseteq\overline{S}\end{array}\end{equation}and
            \begin{equation}\label{RarS}\begin{array}{l}
            \overline{Ra}_{r}\subseteq\overline{S}\Rightarrow
            (W_{r})^{-1}\overline{Ra}_{r}=\overline{Ra}^{r}
            \subseteq\overline{S}\\\hspace{1cm}\Rightarrow
            (W^{r})^{-1}\overline{Ra}^{r}=\overline{Ra}
            \subseteq\overline{S}.\end{array}\end{equation} This
            proves the theorem $\blacksquare.$

            For rational number terms, Theorem \ref{Sterms} holds
            here.  From Eq. \ref{termpxx} one sees that for rational
            number structures, \begin{equation}\label{termRa}
            t^{r}=(\sum_{j,k=1}^{m})^{r}\frac{(ra)^{j}}
            {(rb)^{k}}\mbox{}^{r}=r(\sum_{j,k=1}^{m})
            \frac{(a)^{j}}{(b)^{k}}\mbox{}=rt.\end{equation}

            \section{Real Numbers}\label{RN}
            The description for real numbers is similar to that for
            the rational numbers. The structures $\overline{Ra}$ and
            $\overline{Ra}_{r},$ Eqs. \ref{Ra} and \ref{Rar}, become
            \begin{equation}\label{R}\overline{R}=\{R,+,-,\times,
            \div,<, 0,1\}\end{equation} and \begin{equation}
            \label{Rr}\overline{R}_{r}=\{R,+_{r},-_{r},\times_{r},
            \div_{r},<_{r}, 0_{r},1_{r}\}.\end{equation} The
            external representation of the structure,  whose internal
            representation is $\overline{R}_{r},$  is given in
            terms of the elements, operations, relations and
            constants of $\overline{R}.$ It is denoted by
            $\overline{R}^{r}$ where\begin{equation}
            \label{RrR}\overline{R}^{r}=\{R,+,-,\frac{\times}
            {r}, r\div,<, 0,r\}.\end{equation} Here $r$ is any
            positive real number value in $\overline{R}.$.
            If $r<0$  then $<$ in Eq. \ref{RrR} is replaced by $>.$

            The axioms for the real numbers \cite{real} are similar
            to those for the rational numbers in that both number
            types satisfy the axioms for an ordered field. For this
            reason the proof that $\overline{R}_{r}$ satisfies
            the ordered field axioms if and only if $\overline{R}^{r}$
            satisfies the axioms if and only if $\overline{R}$ does
            will not be given as it is essentially the same as that
            for the rational numbers.

            Real numbers are required to satisfy an axiom of completeness.
            For this axiom, let $\{(a_{r})_{j}\}=\{(a_{r})_{j}:
            j=0,1,2,\cdots\}$ be a sequence of real numbers in
            $\overline{R}_{r},$ Eq. \ref{Rr}. Let  $r$ be a
            positive real number value in $\overline{R}$.
            $\{(a_{r})_{j}\}$ converges in
            $\overline{R}_{r}$ if\begin{equation}\label{convRr}
            \begin{array}{c}\mbox{For all $\epsilon_{r}>_{r}0$ there
            exists an $h$ such that }\mbox{for all $j,m >h$}\\\hspace{1cm}
            |(a_{r})_{j}-_{r}(a_{r})_{m}|_{r}<_{r}\epsilon_{r}.\end{array}
            \end{equation} Here $|(a_{r})_{j}-_{r}(a_{r})_{m}|_{r}$
            denotes the absolute value, in $\overline{R}_{r},$
            of the difference between $(a_{r})_{j}$ and $(a_{r})_{m}.$

            The numerical value in $\overline{R}^{r},$ Eq. \ref{RrR},
            that is the same as $|(a_{r})_{j}-_{r}(a_{r})_{m}|_{r}$
            is in $\overline{R}_{r},$  is given by $|ra_{j}-ra_{m}|.$
            This is the same value as $|a_{j}-a_{m}|$ is in $\overline{R}.$
            $|ra_{j}-ra_{m}|$ also corresponds to a number value in
            $\overline{R}$ given by\begin{equation}\label{rajram}
            |ra_{j}-ra_{m}|=r|a_{j}-a_{m}|.\end{equation}
            \newtheorem{conv}[Sax]{Theorem}\label{conv}
            \begin{conv}Let $r\neq 0$ be a real number value in
            $\overline{R}.$ The sequence $\{(a_{r})_{j}\}$
            converges in $\overline{R}_{r}$ if and only if
            $\{ra_{j}\}$ converges in $\overline{R}^{r}$ if and only
            if $\{a_{j}\}$ converges in $\overline{R}.$\end{conv}
            \underline{Proof}:\\\hspace{.5cm}The proof is in two parts:
            first $r>0,$ and then $r<0.$

            $r>0:$ Let $\epsilon_{r}$ be a positive number value in
            $\overline{R}_{r}$ such that
            $|(a_{r})_{j}-_{r}(a_{r})_{m}|_{r}<_{r}\epsilon_{r}.$ It follows
            that $|ra_{j}-ra_{m}|<r\epsilon$ in both $\overline{R}^{r}$
             and $\overline{R}.$ Eq. \ref{rajram} gives the result
             that $|a_{j}-a_{m}|<\epsilon$ in $\overline{R}.$

            Conversely Let $|a_{j}-a_{m}|<\epsilon$ be true in
            $\overline{R}.$ Then $|ra_{j}-ra_{m}|<r\epsilon$ is true in
            both $\overline{R}$ and $\overline{R}^{r},$
            and $|(a_{r})_{j}-_{r}(a_{r})_{m}|_{r}<_{r}\epsilon_{r}$ is
            true in $\overline{R}_{r}.$ From this one has the
            equivalences, $$|(a_{r})_{j}-_{r}(a_{r})_{m}|_{r}<_{r}
            \epsilon_{r}\Leftrightarrow r|a_{j}-a_{m}|<r\epsilon\Leftrightarrow
            |a_{j}-a_{m}|<\epsilon.$$ It follows that $\{(a_{r})_{j}\}$
            converges in $\overline{R}_{r}$ if and only if $\{ra_{j}\}$
            converges in $\overline{R}^{r}$  if and only
            if $\{a_{j}\}$ converges in $\overline{R}.$

            $r<0$: It is sufficient to set $r=-1.$ In this case
            $\overline{R}_{-1}$  and $\overline{R}^{-1}$ are
            given by Eqs. \ref{Rr} and \ref{RrR} with $r=-1.$
            In this case\begin{equation}\label{Rm1}\overline{R}_{-1}
            =\{R,+_{-1},-_{-1},\times_{-1},\div_{-1},<_{-1}, 0_{-1}
            ,1_{-1}\}\end{equation}and\begin{equation}
            \label{Rm1R}\overline{R}^{-1}=\{R,+,-,\frac{\times}{-1},
            -1\div,>, 0,-1\}.\end{equation} As was the case for
            integers, and is the case for rational numbers, this
            structure can also be considered as a reflection of
            the structure, $\overline{R},$ through the origin at $0.$

            As before let $\{(a_{-1})_{j}\}$ be a convergent sequence in
            $\overline{R}_{-1},$ Eq. \ref{Rm1}. The statement of convergence is
            given by Eq. \ref{convRr} where $r=-1.$ The statement
            $|(a_{-1})_{j}-_{-1}(a_{-1})_{m}|_{-1}<_{-1}\epsilon_{-1}$
            says that $|(a_{-1})_{j}-_{-1}(a_{-1})_{m}|_{-1}$ is a
            positive number in $\overline{R}_{-1}$ that is less than
            the positive number $\epsilon_{-1}$ and $\geq 0.$

            The corresponding statement in $\overline{R}^{-1},$ Eq.
            \ref{Rm1R}, is\begin{equation}\label{convm1}
            |-a_{j}- (-a_{m})|^{-1}>-\epsilon.\end{equation} Since
            $\overline{R}^{-1}$ is a reflection of
            $\overline{R}$ about the origin, the absolute value, $|-|,$
            in $\overline{R}$ becomes $r|x|=-|x|$ in $\overline{R}^{-1}.$
            In $\overline{R}^{-1},$ the absolute value, $|-|^{-1},$ is always
            positive even though it is always negative in $\overline{R}.$
            For example,   $|-a_{j}-(-a_{m})|^{-1}=-|a_{j}- a_{m}|$
            where $|-|^{-1}=-|-|$ and $|-|$ are the respective absolute
            values in $\overline{R}^{-1}$ and $\overline{R}.$

            In this case Eq. \ref{convm1} can be recast as
            \begin{equation}\label{convm1R}-|a_{j}- a_{m}|>
            -\epsilon.\end{equation} This can be used as the
            convergence condition in Eq. \ref{convRr}. Note that
            $0\geq -|a_{j}- a_{m}|$ and that $\geq$ denotes "less
            than or equal to" in $\overline{R}^{-1}.$

            It follows from this that for $r=-1$ that the sequence
            $\{(a_{-1})_{j}\}$ converges in $\overline{R}_{-1},$ Eq,
            \ref{Rm1}, if and only if the sequence $\{-a_{j}\}$
            converges in $\overline{R}^{-1},$ Eq. \ref{Rm1R}, if and
            only if the sequence $\{a_{j}\}$ converges in
            $\overline{R}.$

            Extension to arbitrary negative $r$ can be done in two
            steps. One first carries out the reflection with $r=-1$.
            This is followed by a scaling with a positive value of
            $r$ as has already been described in Section \ref{GD}.
            $\blacksquare$

            \newtheorem{corr}[Sax]{Theorem}\label{comp}
            \begin{corr}$\overline{R}_{r}$ is complete if and only
            if $\overline{R}^{r}$ is complete if and only if
            $\overline{R}$ is complete.\end{corr}.
            \underline{Proof}:
            Assume that $\overline{R}_{r}$  is complete and that the sequence
            $\{(a_{r})_{j}\}$ converges in $\overline{R}_{r}$. Then
            there is a number value  $\mu_{r}$ in $\overline{R}_{r}$ such that
            $\lim_{j\rightarrow\infty}(a_{r})_{j}=\mu_{r}.$
            The properties of convergence, Eq. \ref{convRr},
            with $(a_{r})_{m}$ replaced by $\mu_{r},$ and Theorem
            \ref{conv}, can be used to show that
            \begin{equation}\label{limRr}
            \lim_{j\rightarrow\infty}(a_{r})_{j}=\mu_{r}\Rightarrow
            \lim_{j\rightarrow\infty}ra_{j}=r\mu\Rightarrow\lim_{j
            \rightarrow\infty}a_{j}=\mu.\end{equation} Here
            $\mu_{r}$ is the same number value in $\overline{R}_{r}$ as
            $r\mu$ is in $\overline{R}^{r}$ as $\mu$ is in
            $\overline{R}.$

            Conversely assume that $\overline{R}$ is complete
            and that $\{a_{j}\}$ converges in $\overline{R}$ to
            a number value $\mu.$ Repeating the
            above argument gives\begin{equation}\label{limLr}
            \lim_{j\rightarrow\infty}(a_{r})_{j}=\mu_{r}\Leftarrow
            \lim_{j\rightarrow\infty}ra_{j}=r\mu\Leftarrow\lim_{j
            \rightarrow\infty}a_{j}=\mu.\end{equation}

            This shows that $\overline{R}_{r}$
            is complete if and only if $\overline{R}^{r}$ is complete if
            and only if $\overline{R}$ is complete. $\blacksquare$

            Since real number structures are fields, Eq. \ref{termpxxtx}
            holds for real number terms. That is \begin{equation}\label{Rterm}
            t^{r}=(\sum_{j,k=1}^{m})^{r}\frac{(ra)^{j}}{(rb)^{k}}\mbox{}^{r}
            =r(\sum_{j,k=1}^{m})\frac{(a)^{j}}{(b)^{k}}\mbox{}=rt.\end{equation}

            These terms can be used to give relations between power series
            in $\overline{R}_{r},\overline{R}^{r},$ and $\overline{R}.$
            Let $P_{r}(n,x_{r})=\sum_{j=1}^{n}(a_{r})_{j}x_{r}^{j}$
            be a power series in $\overline{R}_{r}.$ Then
            $P^{r}(n,rx)$ and $P(n,x)$ are the same power series in
            $\overline{R}^{r}$ and $\overline{R}$ as
            $P_{r}(n,x_{r})$ is in $\overline{R}^{r}.$ This means
            that for each real number value $x_{r},$  $P^{r}(n,rx)$
            and $P(n,x)$ are the respective same number values in
            $\overline{R}^{r}$ and $\overline{R}$ as
            $P_{r}(n,x_{r})$ is in $\overline{R}_{r}.$ Here $rx$ and
            $x$ are the same number values in $\overline{R}^{r}$ and
            $\overline{R}$ as $x_{r}$ is in $\overline{R}_{r}.$

            However, the power series in $\overline{R}$ that corresponds
            to $P_{r}(n,x_{r})$ and $P^{r}(n,rx)$ in
            $\overline{R}_{r}$ and $\overline{R}^{r}$ is obtained
            from Eq. \ref{termpxxtx}. It is \begin{equation}\label{Prnx}
            P_{r}(n,x_{r})=P^{r}(n,rx)=rP(n,x). \end{equation}This
            shows that the element of $R$ that has value
            $P_{r}(n,x_{r})$ in $\overline{R}_{r}$ has value
            $rP(n,x)$ in $\overline{R}.$\footnote{Recall that
            $\overline{R}_{r}$ and $\overline{R}^{r}$ are internal
            and external views of the same structure. The structure
            differs from  $\overline{R}$ by the scaling
            factor, $r.$ Also correspondence is a different concept
            from sameness unless $r=1.$}

            These relations extend to convergent power series.
            Theorem \ref{conv} gives the result that
            $P_{r}(n,x_{r})$ is convergent
            in $\overline{R}_{r}$ if and only if $P^{r}(n,rx)$ is
            convergent in $\overline{R}^{r}$ if and only if $rP(n,x)$
            is  convergent in $\overline{R}.$ It follows from
            Theorem \ref{comp} and Eqs. \ref{limRr} and
            \ref{limLr} that,\begin{equation}\label{limPr}
            \begin{array}{l}\lim_{n\rightarrow\infty}P_{r}(n,x_{r})
            =f_{r}(x_{r})\Leftrightarrow\lim_{n\rightarrow\infty}
            P^{r}(n,rx)=f^{r}(rx)\\\hspace{1cm}\Leftrightarrow
            \lim_{n\rightarrow\infty}rP(n,x)=rf(x).\end{array}
            \end{equation} Here $f_{r}(x_{r})$ is the same analytic
            function \cite{Rudin} in $\overline{R}_{r}$ as $f^{r}(rx)$
            is in $\overline{R}^{r}$ as $f(x)$ is in $\overline{R}.$

             Eqs. \ref{Prnx} and \ref{limPr}  give the result that,
            for any analytic function $f$ , \begin{equation}
            \label{analf} f_{r}(x_{r})=f^{r}(rx)=rf(x).\end{equation} Here
            $f_{r}$ and $f$  are functions in $\overline{R}_{r}$
            and $\overline{R}$ and $f_{r}(x_{r})$ is the same number
            value in $\overline{R}_{r}$ as $f^{r}(rx)=rf(x)$ is in
            $\overline{R}^{r}$ as $f(x)$ is in $\overline{R}.$
            Simple examples are $e_{r}^{x_{r}}=(e^{r})^{rx}=re^{x}$ for the
            exponential and $\sin_{r}(x_{r})=sin^{r}(rx) =r\sin(x)$ for the sine
            function. Caution: $\sin^{2}_{r}(x_{r})=r\sin^{2}(x),$ not
            $r^{2}\sin^{2}(x).$

            \section{Complex Numbers}\label{CN}
            The descriptions of structures for complex numbers is
            similar to that for the real numbers. Let $\overline{C}$
            denote the complex number structure\begin{equation}
            \label{C}\overline{C}=\{C,+,-,\times,\div,^{*},0,1\}.
            \end{equation}For each complex number $c$ let $\overline{C}_{c}$
            be the internal representation of another structure where
            \begin{equation}\label{Cc}\overline{C}_{c}=\{C,+_{c},-_{c},
            \times_{c},\div_{c},^{*_{c}},0_{c},1_{c}\}.\end{equation}

            The external representation of the structure, in terms of
            operations and constants in $\overline{C},$ is given by
            \begin{equation}\label{CcC}\overline{C}^{c}=\{C,+,-,
            \frac{\times}{c},c\div,c(-)^{*},0,c1\}.\end{equation}
            The relations for the field operations are the same as
            those for the real numbers except that $c$ replaces $r.$
            It follows that Eqs. \ref{Rterm} - \ref{analf} hold with $c$
            replacing $r.$ These equations show that analytic functions,
            $f_{c}(x_{c})$ in $\overline{C}_{c},$ have corresponding
            functions in $\overline{C}$ given by\begin{equation}\label{fcxc}
            f_{c}(x_{c})=cf(x).\end{equation} Here $x_{c}$ denotes
            the same number in $\overline{C}_{c}$ as $x$ is in
            $\overline{C}.$

            The relation for complex conjugation is given
            by $a_{c}^{*_{c}}=ca^{*}$ It is \emph{not} $a_{c}^{*_{c}}
            =c^{*}a^{*}.$  One way to show this is through the
            requirement that the relation for complex conjugation
            must be such that $1_{c}$ is a real number value in
            $\overline{C}_{c}$ if and only if $1$ is a real number
            value in $\overline{C}.$ This requires that the
            equivalences $$1_{c}^{*_{c}}=1_{c}\Leftrightarrow
            (c1)^{*_{c}}=c1\Leftrightarrow c(1^{*})=c1\Leftrightarrow
            1^{*}=1$$ be satisfied. These equivalences show that
            $(c1)^{*_{c}}=c(1^{*})$ or more generally\begin{equation}
            \label{acomconj}(a_{c})^{*_{c}}=(ca)^{*_{c}}=c(a^{*}).
            \end{equation} Note that any value for $a$ is possible
            including $c$ or powers of $c.$ For example,
            $(c_{c}^{n})^{*_{c}}=(c(c^{n}))^{*_{c}}=c(c^{n})^{*}.$

            As values of elements of the base set, $C$, Eq.
            \ref{acomconj} shows that the element of $C$ that has
            value $a_{c}^{*_{c}}$ in $\overline{C}_{c}$ has value $ca^{*}$
            in $\overline{C}.$ This is different from the element of
            $C$ that has the same value, $a^{*},$ in $\overline{C}$
            as $a_{c}^{*_{c}}$ is in $\overline{C}_{c}.$

            Another representation of the relation of complex
            conjugation in $\overline{C}_{c}$ to that in
            $\overline{C}$ is obtained by writing $c=|c|e^{i\phi}.$
            Here $|c|$ is the absolute value of $c.$ This can be
            used to write\begin{equation}\label{acomconjphi}
            (a_{c})^{*_{c}}=(ca)^{*_{c}}=e^{2i\phi}c^{*}a^{*}.
            \end{equation} That is
            $(c-)^{*_{c}}=e^{2i\phi}c^{*}(-)^{*}.$

            For most of the axioms, proofs that $\overline{C}_{c}$
            satisfies an axiom if and only if $\overline{C}$ does
            are similar to those for the number types already
            treated.   However, it is worth discussing some of the
            new axioms.  For example the complex conjugation axiom
            \cite{comcon} $(x^{*})^{*}=x$ has an easy proof.  It is
            based on Eq. \ref{acomconj}, which gives $$(a_{c}^{*_{c}})^{*_{c}}
            =((ca)^{*_{c}})^{*_{c}}=(c(a^{*}))^{*_{c}}=c(a^{*})^{*}.$$
            From this one has the equivalences,
            $$(a_{c}^{*_{c}})^{*_{c}}=a_{c}\Leftrightarrow
            c(a^{*})^{*}=ca\Leftrightarrow (a^{*})^{*}=a.$$ This
            shows that $(x^{*})^{*}=x$ is valid for
            $\overline{C}_{c}$ if and only if it is valid for
            $\overline{C}^{c}$ and $\overline{C}.$

            The other axiom to consider is that for algebraic
            closure.\newtheorem{algclos}[Sax]{Theorem}
            \begin{algclos}$\overline{C}_{c}$ is algebraically
            closed if and only if $\overline{C}^{c}$ is algebraically closed
            if and only if $\overline{C}$ is algebraically closed.
            \end{algclos}.\underline{Proof}:\\ The proof consists in
            showing that any polynomial equation has a solution in
            $\overline{C}_{c}$ if and only if the same polynomial
            equations in $\overline{C}^{c}$ and $\overline{C}$ have
            the same solutions.

            Let $\sum_{j=0}^{n}b_{j}x^{j}=0$ denote a polynomial
            equation that has a solution $a_{c}$ in $\overline{C}_{c}.$
            Then $\sum_{j=0}^{n}(b_{c})_{j}(a_{c})^{j}=0$ in
            $\overline{C}_{c}.$ Carrying out the replacements
            $(b_{c})_{j}=cb_{j},$ $a_{c}=ca,$  $\times_{c}=\times/c,$ and
            $+_{c}=+$ gives the implications $$\begin{array}{l}
            (\sum_{j=0}^{n})_{c}(b_{c})_{j}(a_{c})^{j}=0\Rightarrow
            (\sum_{j=0}^{n})^{c}(cb_{j}(ca)^{j})^{c}=0\\\hspace{1cm}\Rightarrow
            \sum_{j=0}^{n} cb_{j}a^{j}=0\Rightarrow c\sum_{j=0}^{n}b_{j}
            a^{j}=0\Rightarrow\sum_{j=0}^{n} b_{j}a^{j}=0.\end{array}$$
            From the left, the first equation is in $\overline{C}_{c},$
            the second  is in $\overline{C}^{c}$, and the other three
            are in $\overline{C}.$ This shows that if $a_{c}$
            is the solution of a polynomial in $\overline{C}_{c}$
            then $ca$ is the solution of the same polynomial equation in
            $\overline{C}^{c},$ and  $a$ is the solution of the same
            polynomial equation in $\overline{C}.$

            The proof in the other
            direction consists in assuming that $\sum_{j=0}^{n}b_{j}
            a^{j}=0$  is valid in $\overline{C}$  for some number
            $a,$ and reversing the implication directions  to obtain
            $(\sum_{j=0}^{n})_{c}(b_{c})_{j}(a_{c})^{j}=0$ in
            $\overline{C}_{c}.$ $\blacksquare$

            \section{Number Types as Substructures of $\overline{C}_{c},$
            $\overline{C}^{c},$ and $\overline{C}.$}\label{NTSC}
            As is well known each complex number structure contains
            substructures for the real numbers, the rational
            numbers, the integers, and the natural numbers. The
            structures are nested in the sense that each real number
            structure contains substructures for the rational
            numbers, the integers, and the natural numbers, etc. Here
            it is sufficient to limit consideration to the real
            number substructures of $\overline{C}_{c},$ $\overline{C}^{c},$
            and $\overline{C}$.

            To this end let $c$ be a complex number and let $\overline{R}_{c}$
            and $\overline{R}$ be real number substructures in
            $\overline{C}_{c}$ and $\overline{C}.$
            In this case\begin{equation}\label{Rcc}\overline{R}_{c}=
            \{R_{c},+_{c},-_{c},\times_{c},\div_{c},<_{c},0_{c},
            1_{c}\}\end{equation} and \begin{equation}
            \label{Rc}\overline{R}=\{R,+,-,\times,\div,<,0,1\}.
            \end{equation}

            The field operations in $\overline{R}_{c}$ and
            $\overline{R}$ are  the same as those in
            $\overline{C}_{c}$ and $\overline{C}$ respectively,
            restricted to the number values of the base sets, $R_{c}$
            and $R.$

            Expression of the field operations and constants of
            $\overline{R}_{c}$ in terms of those of $\overline{C}$
            gives a representation of $\overline{R}_{c}$ that
            corresponds to Eq. \ref{CcC}. It is\begin{equation}
            \label{RccR}\overline{R}^{c}=\{R_{c},+,-,\frac{\times}{c},
            c\div,<^{c},0,c1\}.\end{equation}

            Recall that the number values associated with the
            elements  of a base set are not fixed but depend on the structure
            containing them. $R$ contains just those elements of $C$
            that have real values in $\overline{C}.$  $R_{c}$ contains
            just those values of $C$ that that have real values in
            $\overline{C}^{c}$ and $\overline{C}_{c}.$

            Let $x$ be an element of $R$ that has a real value, $a,$ in
            $\overline{R}.$ This is different from another element, $y,$
            of $R_{c}$ that has the same real value, $a_{c},$ in
            $\overline{R}_{c}$ as $a$ is in $\overline{R}.$ The element,
            $y,$ also has the value, $ca,$ in $\overline{R}^{c},$ which is
            the same real value in $\overline{R}^{c}$ as $a$ is in $\overline{R}.$

            This shows that $y$ cannot be an element of $R$ if $c$ is complex.
            The reason is that $ca$ is a complex number value in $\overline{C},$
            and $R$ cannot contain any elements of $C$ that have complex values
            in $\overline{C}.$

            It follows that, if $c$ is complex, then $R_{c}$ and $R$
            have no elements in common, except for the element with
            value $0$, which is the same in all structures.

            The order relations $<_{c},$  $<^{c},$ and $<$ are defined
            relations as ordering is not a basic relation for
            complex numbers. A simple definition of $<$ in $\overline{R}$ is
            provided by defining $a<b$ in $\overline{C}$ by\begin{equation}
            \label{LT}a<b\left\{\parbox{6cm}{if $a$
            and $b$ are real number values in $\overline{C}$ and
            there exists a positive real number value, $d,$ in
            $\overline{C}$ such that $a+d=b.$}\right.\end{equation} Here
            $d$ is a positive real number in $\overline{C}$
            if there exists a number  $g$ in $\overline{C}$
            such that $d=(g^{*})g.$ The definition for $<_{c}$ in
            $\overline{R}_{c}$ is obtained by putting $c$
            subscripts everywhere in Eq. \ref{LT}.

            The order relation $<^{c}$  in $\overline{R}^{c}$ is
            \begin{equation}\label{LTc}ca<^{c}cb\left\{
            \parbox{6cm}{if $ca$ and $cb$
            are real numbers in $\overline{C}^{c},$ and there exists a
            positive real number $cd$ such that $ca+
            cd=cb.$}\right.\end{equation}Here $cd$ is a
            positive real number in $\overline{C}^{c}$ if there exists
            a number  $cg$ in $\overline{C}^{c}$ such that
            $cd=(cg)^{*_{c}}cg.$

            One still has to prove that these definitions of ordering
            are equivalent:
            \newtheorem{order}[Sax]{Theorem}\begin{order}Let $a,$ $b$
            be real number values in $\overline{R}.$ Then $a<b
            \Leftrightarrow ca<^{c}cb\Leftrightarrow a_{c}<_{c}b_{c}.$
            \end{order}\underline{Proof}:\\ Replace the three order
            statements by their definitions, Eqs. \ref{LT} and \ref{LTc},
            to get $a+d=b\Leftrightarrow ca+cd=cb\Leftrightarrow a_{c}+_{c}
            d_{c}=b_{c}.$ Here $d=g^{*}g$, $cd=(cg)^{*^{c}}(cg)$, and
            $d_{c}=g_{c}^{*_{c}}g_{c}.$

            Let $a,b,g$ be real numbers\ values in $\overline{R}.$
            Then $ca,cb,cg$ and $a_{c},b_{c},g_{c}$ are the
            same real number values in $\overline{R}^{c}$ and
            $\overline{R}_{c}$ as $a,b,g$ are in $\overline{R}.$

            Define the positive number value $d$ by $d=g^{*}g.$ Then by
            Eq. \ref{CcC}, $g^{*}g$ is the same number value in
            $\overline{R}$ as $$(cg)^{*_{c}}\frac{\times}{c}(cg)=
            cg^{*}g=cd$$ is in $\overline{R}^{c}$ as $g_{c}^{*_{c}}
            g_{c}=d_{c}$ is in $\overline{R}_{c}.$ Thus $cd$ and $d_{c}$
            are the same positive number values in $\overline{R}^{c}$
            and $\overline{R}_{c}$ as $d$ is in $\overline{R}.$

            From this one has $$a+d=b\Rightarrow ca+cd=cb\Rightarrow
            a_{c}+_{c} d_{c}=b_{c},$$ which gives $$a<b\Rightarrow
            ca<^{c}cb\Rightarrow a_{c}<_{c}b_{c}.$$ The reverse implications
            are proved by starting with $a_{c},b_{c}$, and
            $d_{c}=g_{c}^{*_{c}}g_{c}$  and using an argument
            similar to that given above. $\blacksquare$

            \section{Discussion}\label{D} So far the existence of
            many different isomorphic structure representations of
            each number type has been shown. Some properties of
            the different representations, such as  the fact that
            number values, operations, and relations, in one
            representation are related to values in another
            representation by scale factors have been described.

            Here some additional aspects of these structure
            representations  and their effect on other areas
            of mathematics will be briefly summarized.

            One aspect worth noting is the fact that the
            rational, real, and complex numbers are a multiplicative
            group.  This can be used to define multiplicative
            operations on the structures of these three number types
            and use them to define groups whose elements are the scaled
            structures.

            These are based on  maps from the number values to
            scaled representations.  The maps are listed below for each
            type of number:
            \begin{itemize}
            \item $r\rightarrow \overline{Ra}_{r} \mbox{ $r$ a rational number}$
            \item $r\rightarrow \overline{R}_{r} \mbox{ $r$ a real number}$
            \item $c\rightarrow \overline{C}_{c} \mbox{ $c$ a complex number}.$
            \end{itemize} These maps are restricted to nonzero values\footnote{It
            may be useful in some cases to include empty structures associated
            with the number $0.$ These are represented by the map
            extensions $0\rightarrow \overline{S}_{0}$ where
            $S=Ra,R,C.$} of $r,r,c.$

            The group properties of the  fields for rational, real,
            and complex numbers, induce corresponding properties
            in the collection of structures for each of these
            number types.  The discussion will be limited to complex
            numbers since extension to other number types is similar.

            Let $G_{C}$ denote the collection of complex number
            structures that differ by arbitrary nonzero complex
            scaling factors. Define the operation, $\diamond$, by
            \begin{equation}\label{GrpCcd}\overline{C}_{c}\diamond
            \overline{C}_{d}=\overline{C}_{cd}.\end{equation}

            Justification for this definition is provided by the
            description of iteration  that includes Eqs.
            \ref{Sqpdown}-\ref{SqpS}. $G(C)$ is defined relative to
            the identity group structure, $\overline{C},$ as
            $c,d$ and $cd$ are values in $\overline{C}.$ Every
            structure $\overline{C}_{c}$ has an inverse,
            $\overline{C}_{1/c}$ as \begin{equation}\label{CcC1c}
            \overline{C}_{c}\diamond\overline{C}_{1/c}=
            \overline{C}.\end{equation}

            This shows that $G(C)$ is a group of complex number
            structures induced by the multiplicative group of
            complex number values in $\overline{C}.$ Its properties
            mirror the properties of the group of number values.

            Another consequence of the existence of representations
            of number types that differ by scale factors is that
            they affect other mathematical systems that are based
            on different number types as scalars. Examples include any
            system type that is closed under multiplication by
            numerical scalars.   Specific examples are
            vector spaces based on either real or complex scalars,
            and operator algebras.

            For vector spaces based on
            complex scalars, one has for each complex number $c$ a
            corresponding pair of structures, $\overline{V}_{c},
            \overline{C}_{c}.$ The scalars for $\overline{V}_{c}$
            are number values in $\overline{C}_{c}.$

            The representation of $\overline{V}_{c}$ as a
            structure with the basic operations expressed in
            terms of those in $\overline{V},$ with $\overline{C}$
            as scalars, is similar to that for $\overline{C}_{c},$
            Eq. \ref{CcC}. For Hilbert spaces, the structure
            $\overline{H},$ based on $\overline{C},$ is given by
            \begin{equation}\label{H}\overline{H}=\{H,+,
           -,\cdot, \langle -,-\rangle,\psi\}.\end{equation} Here
           $\cdot$ and $\langle-,-\rangle$ denote scalar vector
           multiplication and scalar product. $\psi$ denotes a
           general state in $\overline{H}.$

           The representation of another Hilbert space structure
           that is based on $\overline {C}_{c},$ is given by
           \begin{equation}\label{Hc}\bar{H}_{c}=\{H,
           +_{c},-_{c},\cdot_{c},\langle-,-\rangle_{c},
           \psi_{c}\}.\end{equation} The representation that expresses
           the basic operations of $\overline{H}_{c}$ in terms of
           those for $\overline{H}$ is given by\begin{equation}
           \label{HcH}\overline{H}^{c}=\{H,+,-,\frac{\textstyle
           \cdot}{\textstyle c},\frac{\textstyle\langle -, -\rangle}
           {\textstyle c},c\psi\}.\end{equation} Details on the
           derivation of this for the case in which $c$ is a real
           number are given in \cite{BenNGF}.

           A possible use of these structures in physics is based on
           an approach to gauge theories \cite{Mack,Montvay} in
           which a finite dimensional vector space $\overline{V}_{x}$ is
           associated with each space time point $x.$ So far just
           one complex number field, $\overline{C}$, serves as the
           scalars for all the $\overline{V}_{x}.$

           The possibility of generalizing this approach by
           replacing $\overline{C}$ with different complex number
           fields, $\overline{C}_{x}$ at each point $x$ has been
           explored in \cite{BenNGF}. If \begin{equation}
           \label{Cx}\bar{C}_{x}=\{C_{x},+_{x},-_{x},
           \times_{x},\div_{x},\mbox{}^{*_{x}},0_{x},
           1_{x}\}.\end{equation} and $y=x+\hat{\nu}dx$ is a
           neighbor point of $x,$ then the local representation of
           $\overline{C}_{y}$ at $x$ is given by \begin{equation}
           \label{Cxr}\bar{C}_{r,x}=\{C_{x},+_{r,x},-_{r,x},
           \times_{r,x},\div_{r,x},\mbox{}^{*_{r,x}},0_{r,x},
           1_{r,x}\}.\end{equation} Number values in
           $\overline{C}_{r,x}$ are denoted by $a_{r.x}.$

           $\overline{C}_{r,x}$  corresponds to $\overline{C}_{c}$,
           Eq. \ref{Cc} with $c$ restricted to be a real number.
           The use of $r$ instead of $c$ was done in \cite{BenNGF}
           to keep the presentation as simple as possible.  However
           it not necessary. Here $r=r_{y,x}$ is a space time dependent
           real number associated with the link from $x$ to $y.$

           The local representation of $\overline{C}_{y}$ on
           $\overline{C}_{x},$ corresponds to Eq. \ref{CcC}.
           It is given by \begin{equation}
           \label{CxrCx}\overline{C}^{r}_{x}=\{C_{x},+_{x},-_{x},
           \frac{\times_{x}}{r},r\div_{x},r(-)^{*_{x}},0_{x},
           r1_{x}\}.\end{equation}

           At point, $x,$ $\overline{C}_{x}$ is the complex
           number structure that is available to an observer,
           $O_{x},$ at $x.$ Also $\overline{C}_{y}$ (given by
           Eq. \ref{Cx} with $y$ replacing the subscript $x$)
           is the complex number structure  available to
           $O_{y},$ at point $y.$ The number value, $a_{y},$ is
           the same value for $O_{y}$ as $a_{x}$ is for $O_{x}.$

           However, the value $a_{y}$ at $y$ is seen by $O_{x}$ as
           an element of $\overline{C}^{r}_{x},$ Eq. \ref{CxrCx}.
           This means that $O_{x}$ sees $a_{y}$ as the value
           $ra_{x}$ in $\overline{C}_{x}.$ Another way to express
           this is to refer to $ra_{x}$ as the local representation,
           or correspondent, of $a_{y}$ at $x.$ Then $\overline{C}^{r}_{x}$ is the
           local representation of $\overline{C}_{y}$ on
           $\overline{C}_{x}.$ Recall that $ra_{x}=a_{r,x}$ is the same
           number value in $\overline{C}^{r}_{x}$ and $\overline{C}_{r,x}$
           as $a_{x}$ is in $\overline{C}_{x}.$

           It is proposed in \cite{BenNGF} that this space time dependence
           of complex number structures may play a role
           in any theory where a comparison is needed of values of a complex
           valued function or field $f(x)$ at different space time points.
           An example is the space time derivative in direction $\hat{\mu}$
           of $f:$\begin{equation}\label{dmux}\partial_{\mu,x}f=
           \frac{f(x+dx^{\mu})-f(x)}{\partial x^{\mu}}.\end{equation}

           The subtraction in Eq. \ref{dmux} has meaning if and only if
           both terms are in the same complex number structure, such as
           $\overline{C}_{x}.$ This is achieved by defining a covariant
           derivative\begin{equation}\label{Dmux}D_{\mu,x}f=
           \frac{r_{\mu,x}f(x+dx^{\mu})_{x}-f(x)}{\partial
           x^{\mu}}.\end{equation} Here $r_{\mu,x}f(x+dx^{\mu})_{x},$
           as a number value in $\overline{C}_{x},$ is the local
           representation or corresponding value of $f(x+dx^{\mu})$ in
           $\overline{C}_{x+dx^{\mu}}.$ Also  $f(x+dx^{\mu})_{x}$ is the
           same number value in $\overline{C}_{x}$ as $f(x+dx^{\mu})$
           is in $\overline{C}_{x+dx^{\mu}}.$

           Gauge fields are introduced by expressing the real number
           $r_{y,x}$ by\begin{equation}\label{ryxAx}r_{y,x}=
           e^{\vec{A}(x)\cdot\hat{\nu}dx}=e^{\sum_{\mu}A_{\mu}(x)
           dx^{\mu}}=\prod_{\mu}r_{\mu,x}.\end{equation} Use of Eqs.
           \ref{Dmux} and \ref{ryxAx} in gauge theory Lagrangians,
           expanding exponentials to first order in small terms,
           and keeping terms that are invariant under local gauge
           transformations, results in $\vec{A}(x)$ appearing as a
           gauge boson that can have mass (i.e. a mass term is
           optional in the Lagrangians.)

           At present there are no immediate and obvious uses of
           $\vec{A}(x)$ as a gauge field in physics. It is suspected
           that $\vec{A}(x)$ may be relevant to some way to relate
           mathematics and physics at a basic level. As such it
           might be useful for development of a coherent theory
           of physics and mathematics together \cite{BenCTPM,BenCTPM1}.
           More work needs to be done to see if these ideas
           have any merit.

             \section*{Acknowledgement}
          This work was supported by the U.S. Department of Energy,
          Office of Nuclear Physics, under Contract No.
          DE-AC02-06CH11357.

            \end{document}